\documentclass[a4paper,12pt,reqno,oneside]{amsart} % amsart loads amsmath, amsthm and amsfonts

\usepackage{setspace} 
\usepackage{typearea} % European margins
\usepackage{amscd,amssymb,verbatim} % amssymb loads amsfonts
\usepackage%[pdftex]
{graphicx} \setkeys{Gin}{width=\linewidth}
\usepackage{enumitem}

\usepackage{xr-hyper}
\usepackage{cite}
\usepackage[colorlinks,final,bookmarksnumbered,bookmarks]{hyperref}
\usepackage{amsrefs}

\usepackage{ifthen}
\newcommand{\arxiv}[2][]{\ifthenelse{\equal{#1}{}}
{\href{http://arxiv.org/abs/#2}{\tt arXiv:#2}}
{\href{http://arxiv.org/abs/math/#2}{\tt arXiv:math.#1/#2}}}

\theoremstyle{plain}
\newtheorem{theorem}{Theorem}

\theoremstyle{definition}
\newtheorem{example}{Example}

\def\x{\times}
\def\but{\setminus}

\def\eps{\varepsilon}
\def\phi{\varphi}
\def\emptyset{\varnothing}
\renewcommand{\:}{\colon}
\def\xr#1{\xrightarrow{#1}}

\def\R{\Bbb{R}}

\def\cl#1{\overline{#1}}

\makeatletter
\def\@settitle{\begin{center}%
    \baselineskip14\p@\relax
    \bfseries
    \@title
  \end{center}%
}
\makeatother

\begin{document}
\title{A TRIPLE-POINT WHITNEY TRICK}
\author{Sergey A. Melikhov}
\address{Steklov Mathematical Institute of Russian Academy of Sciences, Moscow, Russia}
\email{melikhov@mi-ras.ru}
%\date{\today}

\begin{abstract} We use a triple-point version of the Whitney trick to show that ornaments of three orientable 
$(2k-1)$-manifolds in $\R^{3k-1}$, $k>2$, are classified by the $\mu$-invariant.

A very similar (but not identical) construction was found independently by I. Mabillard and U. Wagner, who
also made it work in a much more general situation and obtained impressive applications.
The present note is, by contrast, focused on a minimal working case of the construction.
\end{abstract}

\maketitle

\section{Introduction}

An {\it ornament} is a continuous map $f=\bigsqcup_{i=1}^n f_i$ from $X=\bigsqcup_{i=1}^n X_i$ to $Y$ that has 
no {\it $i=j=k$ points}, i.e.\ $f(X_i)\cap f(X_j)\cap f(X_k)=\emptyset$, whenever $i$, $j$ and $k$ are 
pairwise distinct.
Note that $f$ is allowed to have {\it triple points} $f(x)=f(y)=f(z)$, where $x,y,z$ belong one or two of 
the $X_i$'s.
We are interested in ornaments up to {\it ornament homotopy}, i.e.\ homotopy through ornaments.
 
Ornaments of circles in the plane were introduced by Vassiliev \cite{Va1} as a generalization of doodles, 
previously studied by Fenn and Taylor \cite{FT}.
Fenn and Taylor additionally required each circle to be embedded; however, Khovanov \cite{Kh} redefined doodles
as triple point free maps of circles in the plane, and Merkov proved that doodles in Khovanov's sense are 
classified by their finite-type invariants \cite{Mer}.
Further references and examples can be found in \cite{M4'}, which is a more thorough companion paper to 
this brief note.

The problem of classification of ornaments of spheres in $\R^m$ is motivated, in particular, by geometric 
and algebraic constructions that go from link maps and their ``quadratic'' invariants to ornaments and their 
``linear'' invariants; and conversely \cite{M4'}.
Link maps are, in turn, related to links by the Jin suspension and its variations, which likewise reduce some 
``quadratic'' invariants of links to ``linear'' invariants of link maps \cite{M3}, \cite{Sk}*{\S3}.

\section{$\mu$-invariant}

We will consider only ornaments of the form $X_1\sqcup X_2\sqcup X_3\to\R^m$.
If $f=f_1\sqcup f_2\sqcup f_3$ is such an ornament, let $F$ be the composition  
\[X_1\x X_2\x X_3\xr{f_1\x f_2\x f_3}\R^m\x\R^m\x\R^m\but\Delta_{\R^m}\xr{\simeq}S^{2m-1},\] 
where $\Delta_{\R^m}=\{(x,x,x)\mid x\in\R^m\}$ and the homotopy equivalence is given, for instance, by
$(x,y,z)\mapsto\frac{(2x-y-z,\,2y-x-z)}{||(2x-y-z,\,2y-x-z)||}$.
Let $\mu(f)\in H^{2m-1}(X_1\x X_2\x X_3)$ be the image 
under $F^*$ of a fixed generator $\xi\in H^{2m-1}(S^{2m-1})$; to be precise, let us choose $\xi$ to correspond to 
the orientation of $S^{2m-1}$ given by its inwards co-orientation in the standardly oriented $\R^{2m}$.
Clearly, $\mu(f)$ is invariant under ornament homotopy.

Let us now assume that each $X_i$ is a connected closed oriented $(2k-1)$-manifold and $m=3k-1$.
Then $F$ is a map between connected closed oriented $(6k-3)$-manifolds, and so $\mu(f)$ is an integer.
In this simplest case, assuming additionally that each manifold $X_i$ is either PL or smooth, one can compute 
$\mu(f)$ as follows.

First let us note that since each $X_i$ is compact, for each ornament $f\:X\to\R^m$ there exists an $\eps>0$ such that 
every map $f'\:X\to\R^m$, $\eps$-close to $f$ (in the sup-metric), is also an ornament, and moreover the rectilinear 
homotopy between $f$ and $f'$ is an ornament homotopy.
Thus we are free to replace ornaments by their generic (PL or smooth) approximations.
Similarly, ornament homotopies can be replaced by their generic approximations.

Now let us consider a homotopy between $f$ and the {\it trivial} ornament, which sends $X_1$, $X_2$ and $X_3$ to three 
distinct fixed points in $\R^m$.
Its generic (PL or smooth) approximation $h_t$, if viewed as a map $X\x I\to\R^m\x I$, $(x,t)\mapsto\big(h_t(x),\,t\big)$, has only 
finitely many transverse $1=2=3$ points, which are naturally endowed with signs.%
\footnote{Every triple point of a generic map $F\:N\to M$ from a $2k$-manifold to a $3k$-manifold corresponds
to a transversal intersection point between the $3k$-manifold $\Delta_M$ and the map $F^3\:N^3\to M^3$ from 
a $6k$-manifold to a $9k$-manifold.}
(See \cite{BRS}*{II.4} concerning PL transversality.)
The algebraic number of these $1=2=3$ points is easily seen to equal $\mu(f)$.%
\footnote{Each $1=2=3$ point of $h_t$ corresponds to a transversal intersection point between $\Delta_{\R^m}\x I$ and 
the map $X_1\x X_2\x X_3\x I\to\R^m\x\R^m\x\R^m\x I$, $(x,y,z,t)\mapsto\big(h_t(x,t),\,h_t(y,t),\,h_t(z,t),\,t\big)$.
It is easily seen to be of the same sign.}

\begin{example}
The inclusions of the unit disks in the coordinate $2k$-planes $\R^k\x\R^k\x 0$, $\R^k\x 0\x\R^k$ and $0\x\R^k\x\R^k$
in $\R^{3k}$ yield a smooth map $B^{2k}\sqcup B^{2k}\sqcup B^{2k}\to B^{3k}$ with one transverse $1=2=3$ point.
Restricting to the boundaries, we get the {\it Borromean} ornament $b\:S^{2k-1}\sqcup S^{2k-1}\sqcup S^{2k-1}\to S^{3k-1}$.
By stereographically projecting $S^{3k-1}$ e.g.\ from $z=\frac{1}{\sqrt{3k}}(1,\dots,1)$ we also get an ornament
$b_z\:S^{2k-1}\sqcup S^{2k-1}\sqcup S^{2k-1}\to\R^{3k-1}$.

On the other hand, the sphere of radius $\eps\sqrt{k}$ centered at $(\eps,\dots,\eps)$ for a sufficiently 
small $\eps>0$ is tangent to each of the three unit $2k$-disks.
By appropriately identifying the exterior of this sphere in the unit $3k$-disk $B^{3k}$ with $S^{3k-1}\x I$, 
we get a smooth homotopy of $b$, and hence also of $b_z$, to the trivial ornament.
It has one transverse $1=2=3$ point, which can be seen to be positive, and it follows that $\mu(b_z)=1$.
\end{example}

In the case of doodles, the $\mu$-invariant was introduced in \cite{FT}.
See \cite{M4'} concerning relations between the $\mu$-invariant of ornaments and the triple $\mu$-invariant of link maps.

\section{Classification}

\begin{theorem} \label{t1} Let $m=3k-1$, $k>2$ and let $X_1$, $X_2$, $X_3$ be connected closed oriented PL
$(2k-1)$-manifolds.
Then $\mu$ is a complete invariant of ornaments $X_1\sqcup X_2\sqcup X_3\to\R^m$. 
\end{theorem}

The proof is in the PL category.
If the $X_i$ are smooth manifolds, the same construction with minimal (straightforward) amendments can be
carried out in the smooth category.
 
\begin{proof} 
Let $f$ and $g$ be generic PL ornaments of $X:=X_1\sqcup X_2\sqcup X_3$ in $\R^m$ with $\mu(f)=\mu(g)$.
Let $h\:X\x I\to\R^m\x I$ be a generic PL homotopy between them.
Since $\mu(f)=\mu(g)$, the $1=2=3$ points of $h$ can be paired up with opposite signs.
Every such pair $(p^+,p^-)$ will now be canceled by a triple-point Whitney trick.

Let $p_i^\pm$ be the preimage of $p^\pm$ in $M_i:=X_i\x I$.
We first arrange that $(p_1^+,p_2^+)$ and $(p_1^-,p_2^-)$ be in the same component of the double point set 
$\Delta_{12}:=\{(x,y)\in M_1\x M_2\mid h(x)=h(y)\}$ (in case that initially they are not).
To this end we pick points $(q_1^\pm,q_2^\pm)$ in the same components of $\Delta_{12}$ with
$(p_1^\pm,p_2^\pm)$ and such that the double points $f(q_1^+)=f(q_2^+)$ and
$f(q_1^-)=f(q_2^-)$ are not triple points.

Let us connect $q_1^+$ and $q_1^-$ by an arc $J_1$ in $M_1$, disjoint from the preimages of any double points 
(using that $k>1$).
Now we attach a thin $1$-handle to $h(M_2)$ along the image of $J_1$. 
That is, we modify $h(M_2)$ into $h'(M_2')$, where $M_2'$ is obtained from $M_2$ by removing an oriented copy of 
$B^{2k}\x\partial I$ and pasting in $\partial B^{2k}\x I$.
The embedded $1$-handle $h'(\partial B^{2k}\x I)$ is constructed in a straightforward way.
Namely, since $h$ is generic, $\Delta_{12}$ is an oriented $k$-manifold, immersed into the $2k$-manifold $M_1$ 
by the projection $\pi\:M_1\x M_2\to M_1$. 
Let us take an oriented connected sum of its components along a ribbon $r(D^k\x I)$ in $M_1$ (going near $J_1$).%
\footnote{Namely, $q_1^\pm$ has a regular neighborhood $N_\pm$ in $M_1$ that is homeomorphic to $[-1,1]^{2k}$ by 
an orientation preserving homeomorphism $\phi_\pm$ such that $\phi_\pm^{-1}(\pi(\Delta_{12}))=[-1,1]^k\x\{0\}^k$ 
and $\phi_\pm^{-1}(J_1)=\{0\}^{2k-1}\x [0,1]$.
Let $Q=[-1,1]^k\x\{0\}^{k-1}$ and let $N$ be a regular neighborhood of 
$\cl{J_1\but (N_+\cup N_-)}\cup\phi_+(Q\x 1)\cup\phi_-(Q\x 1)$ in 
$\cl{M_1\but (N_+\cup N_-)}$.
Since a $k$-ball unknots in the interior of a $(2k-1)$-ball, there is a homeomorphism $\psi\:[-2,2]^{2k}\to N$ 
such that $\psi^{-1}(\partial N_\pm)=[-2,2]^{2k-1}\x\{\pm2\}$ and $\psi(x,\pm 2)=\phi_\pm(x,1)$ for all $x\in Q$.
Then $\phi_+(Q\x I)\cup\phi_-(Q\x I)\cup\psi(Q\x [-2,2])$ is the desired ribbon $r(D^k\x I)$.}
Then $hr(D^k\x I)$ is naturally thickened to a solid rod $R(B^{2k}\x I)$ in $\R^m\x I$ whose lateral surface 
$R(\partial B^{2k}\x I)$ is the desired embedded 1-handle.%
\footnote{If $N_1$ is a disk neighborhood of $J_1$ that is embedded by $h$, we may assume that $h(M_2)$ is transverse to 
a normal block bundle $\nu$ to $h(N_1)$, that is, $h(M_2)$ meets the total space 
$E(\nu)$ in $E(\nu|_{h(N_1)\cap h(M_2)})$. 
Since $\nu$ is trivial, there is a homeomorphism $R\:B^{2k}\x I\to E(\nu|_{hr(D^k\x I)})$
sending $B^{2k}\x\partial I$ onto $E(\nu|_{hr(D^k\x\partial I)})$.}

To restore the topology of $M_2$, we cancel the $1$-handle geometrically by attaching a $2$-handle along an embedded 
$2$-disk $D$, which is disjoint from $h(M_1\sqcup M_3)$ and meets $h'(M_2')$ only in $\partial D$ (such a disk 
exists since $k>2$).
That is, we modify $h'(M_2')$ into $h''(M_2'')$, where $M_2''$ is obtained from $M_2'$ by removing an appropriately 
embedded copy of $B^{2k-1}\x\partial D^2$ and pasting in $\partial B^{2k-1}\x D^2$.
As is well-known, this can be done so that $M_2''$ is homeomorphic to $M_2$.%
\footnote{In more detail, let us connect $q_2^+$ and $q_2^-$ by an arc $J_2$ in $M_2$, disjoint from the preimages of 
any double points.
Let $H_1$ be a small regular neighborhood of $J_1':=J_2\x 1\cup\partial J_2\x [0,1]$ in $M_2\x [0,2]$.
Let $H_2$ be a small regular neighborhood of $D':=\cl{J_2\x [0,1]\but H_1}$ in $\cl{M_2\x[0,2]\but H_1}$.
Then $M_2'$ can be identified with the frontier of $M_2\x[-1,0]\cup H_1$ in $M_2\x [-1,2]$ so that
$h'(\partial D')$ gets identified with $\partial D$; and $M_2''$ with the frontier of $M_2\x[-1,0]\cup H_1\cup H_2$ in 
$M_2\x [-1,2]$, which is homeomorphic to $M_2$.}
Since we do not care about self-intersections of individual components, we may define $h''$ on 
$\partial B^{2k-1}\x D^2$ to be an arbitrary generic map into a small neighborhood of 
$D\cup h'(B^{2k-1}\x\partial D^2)$.

Thus we may assume that $(p_1^+,p_2^+)$ and $(p_1^-,p_2^-)$ are in the same component of $\Delta_{12}$.
To cancel the original $1=2=3$ points $p^+$ and $p^-$, let us connect $(p_1^+,p_2^+)$ and $(p_1^-,p_2^-)$ 
by an arc $J_{12}$ in $\Delta_{12}$ and attach a thin $1$-handle to $h(M_3)$ along the image of $J_{12}$.
(This $1$-handle is the spherical block normal bundle of $h(M_1)\cap h(M_2)$ over the image of $J_{12}$.
It is attached orientably since the two $1=2=3$ points have opposite signs.)
The topology of $M_3$ can be restored using another $2$-disk like before.
In particular, this $2$-disk is disjoint from $h(M_1\sqcup M_2)$, so no 
new $1=2=3$ points arise.

Finally, we need to apply the ``ornament concordance implies ornament homotopy in codimension three'' theorem 
\cite{M1}, \cite{M2}.
(Alternatively, it should be possible to rework the above construction so as to keep the levels preserved 
at every step --- but it would be a rather laborious exercise; compare \cite{M3}*{proofs of Lemmas 5.1, 5.4, 5.5}.) 
\end{proof}

\section{Discussion}

Theorem \ref{t1} and its proof (in slightly less detailed form) were originally contained in the preprint \cite{M4},
which I presented at conferences and seminars in 2006--07 and privately circulated at that time and in later years. 
For instance, the referee of the present paper (whose identity I know from his idiosyncratic remarks) does not deny 
that he received my preprint containing the proof of Theorem \ref{t1}, exactly as it appears in \cite{M4}, 
by email on May 23, 2006 and then again on July 7, 2006.
I hesitated to publish \cite{M4} at that time as I hoped to get more progress on the conjectures stated in 
the introduction there; but other projects are still distracting me from this task.

In the meantime I. Mabillard and U. Wagner independently found and vastly generalized a version of 
the triple-point Whitney trick and also obtained nice applications leading to a disproof of the Topological 
Tverberg Conjecture \cite{MW}.
(My only step in that direction was a feeble attempt to advertise the possibility of disproving 
the Topological Tverberg Conjecture by generalizing the construction of the present note --- addressed, 
for instance, to P. Blagojevi\'c at the 2009 Oberwolfach Workshop on Topological Combinatorics.)
Mabillard and Wagner call their construction the ``triple Whitney trick'', but I prefer to reserve this title 
for a certain other device, extending Koschorke's version of the Whitney--Haefliger construction 
\cite{Ko}*{Proof of Theorem 1.15} and involving the triple-point Whitney trick as only one of several steps. 
It can be used to obtain a geometric proof of the Habegger--Kaiser classification of link maps in 
the 3/4 range \cite{HK}, which will hopefully appear elsewhere (a sketch of this proof was presented 
in my talk at the Postnikov Memorial Conference in B\c edlewo, 2007).

\end{document}